\UseRawInputEncoding 
\documentclass[11pt]{article}

\usepackage{amsmath,amssymb,amsfonts,amsbsy,mathrsfs}

\usepackage[title]{appendix}
\usepackage{titlesec}

\renewcommand{\subsection}{\subsubsection}

\newtheorem{theorem}{Theorem}[section]

\newtheorem{proposition}{Proposition}[section]

\newtheorem{definition}{Definition}[section]

\DeclareMathOperator*{\res}{Res}
\DeclareMathOperator*{\Real}{Re}
\DeclareMathOperator*{\Imag}{Im}

\begin{document}

\title{\bf On weak stability of shock waves \\ in 2D compressible elastodynamics}

\author{{\bf Yuri Trakhinin}\\
Novosibirsk State University, Pirogova str. 1, 630090 Novosibirsk, Russia\\
and\\
Sobolev Institute of Mathematics, Koptyug av. 4, 630090 Novosibirsk, Russia\\
E-mail: trakhin@math.nsc.ru
}

\date{ }

\maketitle

\begin{abstract}
By using an equivalent form of the uniform Lopatinski condition for 1-shocks, we prove that the stability condition found by the energy method in [A. Morando, Y. Trakhinin, P. Trebeschi, Structural stability of shock waves in 2D compressible elastodynamics, {\it Math. Ann.} {\bf 378} (2020) 1471--1504] for the rectilinear shock waves in two-dimensional flows of compressible isentropic inviscid elastic materials is not only sufficient but also necessary for uniform stability (implying structural nonlinear stability of corresponding curved shock waves). The key point of our spectral analysis is a delicate study of the transition between uniform and weak stability. Moreover, we prove that the rectilinear shock waves are never violently unstable, i.e., they are always either uniformly or weakly stable.
\end{abstract}

\section{Introduction}
\label{s1}

We consider the equations of elastodynamics governing the  motion of compressible isentropic inviscid elastic materials. These equations arise as the inviscid limit of the equations of compressible viscoelasticity \cite{Daf,Gurt,Jos} of Oldroyd type \cite{Old1,Old2}. As in \cite{CHW1,CHW2,CHW3,MTT20}, we restrict ourself to two-dimensional (2D) elastic flows and the special case of Hookean linear elasticity. Then, the elastodynamics equations are written as the following system of conservation laws:
\begin{equation}\label{7}
\left\{
\begin{array}{l}
 \partial_t\rho  +{\rm div}\, (\rho v )=0, \\
 \partial_t(\rho v ) +{\rm div}\,(\rho v\otimes v  ) + {\nabla}p-{\rm div}\,(\rho FF^{\top})=0,\\
 \partial_t(\rho F_j ) +{\rm div}\,(\rho F_j\otimes v  -  v\otimes  \rho F_j) =0,\quad j=1,2,
\end{array}
\right.
\end{equation}
where $\rho$ is the density, $v\in\mathbb{R}^2$  is the velocity, $F\in \mathbb{M}(2,2)$ is the deformation gradient,  $F_1=(F_{11},F_{21})$ and $F_2=(F_{12},F_{22})$ are the columns of F, and the pressure $p=p(\rho)$ is a smooth function of $\rho$. System \eqref{7} is supplemented by the identity ${\rm div}\,(\rho F^{\top})=0$ which is the set of the two divergence constraints
\begin{equation}\label{8}
{\rm div}\,(\rho F_j)=0\quad (j=1,2)
\end{equation}
on initial data, i.e., one can show that if equations \eqref{8} are satisfied initially, then they hold for all $t > 0$. Moreover, one can show
\cite[Proposition 1]{QianZhang} that the physical identity (see, e.g., \cite{Daf})
\begin{equation}\label{8'}
\rho\det F =1
\end{equation}
is also a constraint on the initial data for the Cauchy problem. That is, system \eqref{7}, which might seem overdetermined, is actually a closed system for the vector of unknowns  $U=(p,v,F_1,F_2)$.

Taking into account the divergence constraints \eqref{8}, we rewrite \eqref{7} as
\begin{equation}
\left\{
\begin{array}{l}
{\displaystyle\frac{1}{\rho c^2}\,\frac{{\rm d} p}{{\rm d}t} +{\rm div}\,{v} =0,}\\[6pt]
{\displaystyle\rho\, \frac{{\rm d}v}{{\rm d}t}+{\nabla}p -\rho (F_1\cdot\nabla )F_1 -\rho (F_2\cdot\nabla )F_2=0 ,}\\[6pt]
\rho\,{\displaystyle \frac{{\rm d} F_j}{{\rm d}t}-\rho \,(F_j\cdot\nabla )v =0,} \quad j=1,2,
\end{array}\right. \label{9}
\end{equation}
where $c^2=p'(\rho)$ is the square of the sound speed and ${\rm d} /{\rm d} t =\partial_t+({v} \cdot{\nabla} )$ is the material derivative. Equations \eqref{9} form the symmetric system
\begin{equation}
\label{10}
A_0(U )\partial_tU+A_1(U)\partial_1U +A_2(U)\partial_2U=0,
\end{equation}
where $A_0= {\rm diag} (1/(\rho c^2) ,\rho I_6)$,
\[
A_1=\begin{pmatrix}
{\displaystyle\frac{v_1}{\rho c^2}} & e_1 & \underline{0} & \underline{0}  \\[7pt]
e_1^{\top}&\rho v_1I_2 & -\rho F_{11}I_2 & -\rho F_{12}I_2  \\[3pt]
\underline{0}^{\top} &-\rho F_{11}I_2 & \rho v_1I_2 & O_2 \\
\underline{0}^{\top} &-\rho F_{12}I_2 & O_2 & \rho v_1I_2
\end{pmatrix},\quad
A_2=\begin{pmatrix}
{\displaystyle\frac{v_2}{\rho c^2}} & e_2 & \underline{0} & \underline{0}  \\[7pt]
e_2^{\top}&\rho v_2I_2 & -\rho F_{21}I_2 & -\rho F_{22}I_2  \\[3pt]
\underline{0}^{\top} &-\rho F_{21}I_2 & \rho v_2I_2 & O_2 \\
\underline{0}^{\top} &-\rho F_{22}I_2 & O_2 & \rho v_2I_2
\end{pmatrix},
\]
$e_1=(1,0)$, $e_2=(0,1)$, $\underline{0}=(0,0)$ and  $I_m$ and $O_m$ denote the unit and zero matrices of order $m$ respectively. In \eqref{10} we think of the density as a function of the pressure: $\rho =\rho (p)$, $c^2=1/\rho'(p)$. System \eqref{10} is symmetric hyperbolic if   $A_0>0$, i.e.,
\begin{equation}
\rho >0,\quad  \rho '(p)>0. \label{11}
\end{equation}

Let $\Gamma (t)=\{ x_1=\varphi (t,x_2)\}$ be a curve of strong discontinuity for the conservation laws \eqref{7}, i.e., we are interested in solutions of \eqref{7} that are smooth on either side of $\Gamma (t)$. As is known, to be weak solutions such piecewise smooth solutions should satisfy corresponding Rankine-Hugoniot jump conditions at each point of $\Gamma (t)$. The Rankine-Hugoniot conditions for \eqref{7} are written in \cite{CHW1,MTT20}. As in \cite{MTT20}, we are interested in {\it shock waves}. For them there is a non-zero mass transfer flux $\mathfrak{m}$ across the discontinuity curve, and the density has a non-zero jump $[\rho ]$.
Following \cite{MTT20},  from the Rankine-Hugoniot conditions we can deduce the boundary conditions
\begin{equation}\label{bc}
\begin{split}
& [\mathfrak{m}]=0,\quad \mathfrak{M}[V]+\left( 1+(\partial_2\varphi)^2 \right)[p]=0,\quad [v_{\tau}]=0,\\
&  [F_{j\tau}]=0,\quad [\rho F_{j{\rm N}}]=0,\quad j=1,2,
\end{split}
\end{equation}
on the curve $\Gamma (t)$ of a shock wave, where $[g]=g^+|_{\Gamma}-g^-|_{\Gamma}$ denotes the jump of $g$, with $g^{\pm}:=g$ in the domains
\[
\Omega^{\pm}(t)=\{\pm (x_1- \varphi (t,x_2))>0\},
\]
and
\[
\mathfrak{m}^{\pm}=\rho^{\pm} (v_{\rm N}^{\pm}-\partial_t\varphi),\quad \mathfrak{m}=\mathfrak{m}^{\pm}|_{\Gamma},\quad v_{\rm N}^{\pm}=v_1^{\pm}-v_2^{\pm}\partial_2\varphi ,
\]
\[
\mathfrak{M}=\mathfrak{m}^2-(\rho^\pm )^2\left.\left((F_{1{\rm N}}^\pm )^2+(F_{2{\rm N}}^\pm )^2 \right)\right|_{\Gamma},
\quad
F_{j{\rm N}}^{\pm}=F_{1j}^{\pm}-F_{2j}^{\pm}\partial_2\varphi ,
\]
\[
V^\pm =1/\rho^\pm ,\quad
v_{\tau}^{\pm}=v_1^{\pm}\partial_2\varphi +v_2^{\pm},\quad F_{j\tau}^{\pm}=F_{1j}^{\pm}\partial_2\varphi +F_{2j}^{\pm}.
\]
It was assumed in \cite{MTT20} that $\mathfrak{M}\neq 0$ while deducing  \eqref{bc}, but this assumption holds thanks to the Lax conditions \cite{Lax} (see \eqref{25} below).

The free boundary problem for shock waves is the problem for the systems
\begin{equation}
A_0(U^{\pm})\partial_tU^{\pm}+A_1(U^{\pm} )\partial_1U^{\pm}+A_2(U^{\pm} )\partial_2U^{\pm}=0\quad \mbox{in}\ \Omega^{\pm}(t),
\label{21}
\end{equation}
cf. \eqref{10}, with the boundary conditions \eqref{bc} on $\Gamma (t)$ and  the initial
\begin{equation}
{U}^{\pm} (0,{x})={U}_0^{\pm}({x}),\quad {x}\in \Omega^{\pm} (0),\quad \varphi (0,{x}_2)=\varphi _0({x}_2),\quad {x}_2\in\mathbb{R}.\label{indat}
\end{equation}
The initial data \eqref{indat} should satisfy not only the hyperbolicity conditions \eqref{11} but also constraints \eqref{8} and \eqref{8'}. As for the Cauchy problem, one can show that these constraints are preserved by problem \eqref{bc}--\eqref{indat}. That is,
the following proposition holds true.

\begin{proposition}[\cite{MTT20}]
Suppose that problem \eqref{bc}--\eqref{indat} has a smooth solution  $(U^+,U^-,\varphi)$ for $t\in [0,T]$ satisfying the shock wave assumption $\mathfrak{m}\neq 0$. Then, if the initial data \eqref{indat} satisfy \eqref{8} and \eqref{8'} in $\Omega^{\pm} (0)$, then
\begin{equation}\label{8''}
{\rm div}\,(\rho^{\pm} F_j^{\pm})=0\quad \mbox{in}\ \Omega^{\pm}(t) \quad (j=1,2)
\end{equation}
and
\begin{equation}\label{8^}
\rho^{\pm}\det F^{\pm} =1\quad \mbox{in}\ \Omega^{\pm}(t)
\end{equation}
for all $t\in [0,T]$.
\label{p1}
\end{proposition}

The structural (nonlinear) stability of a shock wave means the local-in-time existence and uniqueness of a smooth solution $(U^+,U^-,\varphi)$ to  problem \eqref{bc}--\eqref{indat}. It is well-known that the Lax's $k$--shock conditions \cite{Lax}, which guarantee the correct number of boundary conditions, are necessary for the well-posedness of such free boundary problem. It was shown in \cite{MTT20} that, as in gas dynamics, only {\it  extreme} shocks are possible in elastodynamics. If, without loss of generality, we assume that $\mathfrak{m}>0$, then they are 1-shocks. The 1-shock conditions for \eqref{21} read \cite{MTT20}:
\begin{equation}
\left\{
\begin{array}{l}
v_{\rm N}^--\partial_t\varphi >\sqrt{(\bar{c}^-)^2+(F_{1{\rm N}}^{-})^2+(F_{2{\rm N}}^{-})^2}\,, \\[6pt]
 \sqrt{(F_{1{\rm N}}^{+})^2+(F_{2{\rm N}}^{+})^2}<v_N^+-\partial_t\varphi
 <\sqrt{(\bar{c}^{+})^2+(F_{1{\rm N}}^{+})^2+(F_{2{\rm N}}^{+})^2}  \quad \mbox{on}\ \Gamma (t),
\end{array}
\right.
\label{25}
\end{equation}
where  $\bar{c}^{\pm}= c^{\pm}\sqrt{1+(\partial_2\varphi)^2}$ and $c^{\pm} =1/\sqrt{\rho'(p^{\pm})}$ are the sound speeds ahead and behind of the shock. We see that the Lax conditions \eqref{25} imply the fulfilment of the assumption $\mathfrak{M}\neq 0$.

The study of the structural stability of shock waves in compressible isentropic elastodynamics (for Hookean linear elasticity) was recently begun in \cite{MTT20} for 2D flows. By the energy method based on a symmetrization of the wave equation \cite{BThand,Tsiam} and giving an a priori estimate  without loss of derivatives for solutions of the constant coefficient linearized problem, a condition sufficient for the uniform stability \cite{BS,BThand,M2,Met} of rectilinear shock waves (see \eqref{usc'} below) was found in \cite{MTT20}. As is known, uniform stability, for which we generically have a priori estimates without loss of derivatives,  implies structural stability of corresponding curved shock waves \cite{M2,Met,MZ,Tsiam}.

Moreover, for two particular cases for which the Cauchy--Green stress tensor $\rho {F}{F}^\top$ is diago\-nal for the unperturbed flow, by the test of the uniform Lopatinski condition \cite{Kreiss} it was shown in \cite{MTT20} that the stability condition \eqref{usc'} is necessary for uniform stability. Regarding the case of general deformations, it was also proved in \cite{MTT20} that, as in isentropic gas dynamics \cite{M1}, all compressive shock waves are uniformly stable for convex equations of state (for them the uniform stability condition \eqref{usc'} always holds).

The main goal of the present paper is to show for general deformations that the stability condition \eqref{usc'} found in \cite{MTT20} by the energy method is not only sufficient but also necessary for uniform stability. It is clear that this can be done only by spectral analysis. However, it seems that the classical tools towards the test of the Lopatinski condition and the uniform Lopatinski condition, as it is described in \cite{BS,Kreiss,M1,Met}, are hardly applicable to the model of elastodynamics because of its big technical complexity. We overcome this difficulty by using equivalent definitions of the Lopatinski condition and the uniform Lopatinski condition proposed in \cite{Tcmp} for symmetric hyperbolic systems having the 1-shock property. Roughly speaking, the main idea of \cite{Tcmp} borrowed from \cite{GK} is to carry all the calculations over the unique characteristic incoming in the shock in the region behind of the shock wave.

At the same time, even applying the alternative method from \cite{Tcmp}, we are not able to find an explicit form of the roots of the Lopatinski determinant. Instead of this, motivated by \cite{BRSZ}, we analyze the transition between uniform and weak stability (still using the alternative definition of the uniform Lopatinski condition from \cite{Tcmp}). It turns out that in the space of admissible parameters of the unperturbed flow this transition coincides exactly with the boundary of the parameter domain of uniform stability found in \cite{MTT20} by the energy method. Moreover, we show that the rectilinear shock waves are never violently unstable (this was not proved in \cite{MTT20}). The latter is natural because from the physical point of view the elastic force should play stabilizing role and we know that shock waves in isentropic gas dynamics \cite{M1} are also either uniformly or weakly stable.

The question which is left open in the present paper is that about the deduction of an priori estimate (with a loss of derivatives) for the parameter domain of weak stability. On the one hand, following some constructions of the energy method in \cite{Blmon} for weakly stable shock waves in gas dynamics, we could write down here an energy a priori estimate with a loss of derivatives for the constant coefficient linearized problem, at least, in the particular case when the Cauchy--Green stress tensor $\rho {F}{F}^\top$ is diagonal for the unperturbed flow. However, the peculiarity of the energy method in \cite{Blmon} for weakly stable shocks is such that, roughly speaking, it is unstable against lower-order terms, and so the a priori estimates obtained for the case of constant coefficients cannot be transferred to variable coefficients.

On the other hand, since system \eqref{7} is, in some sense, related to isentropic gas dynamics, there is a natural hope that for it we could get results similar to those in \cite{Col1,Col2}, where energy a priori estimates for weakly stable shock waves in isentropic gas dynamics were derived both for constant and variable coefficients. It is worth noting that the structural stability of these shock waves was then proved in \cite{CS}. Because of the mentioned technical complexity of elastodynamics (even for the 2D case) it is not yet clear whether some assumptions from \cite{Col2} can be checked for it. At the same time, in view of a stabilizing role of the elastic force, it is natural to expect that weakly stable shocks in isentropic elastodynamics are structurally stable. The rigorous proof of this hypothesis is an interesting open problem for a future research.

Regarding characteristic discontinuities in compressible elastodynamics, the linear and structural stability of 2D vortex sheets in isentropic inviscid elastic materials (described by system \eqref{7}) was recently studied in \cite{CHW1,CHW2,CHW3}. We also mention recent results in \cite{CSW} where the structural stability of contact discontinuities in nonisentropic elastodynamics (thermoelasticity), for which the velocity is continuous across the discontinuity surface,  was proved under some stability condition on the piecewise constant background states.

The rest of this paper is organized as follows. In Sec. \ref{s2}, we formulate the constant coefficient linearized  problem and our main result for it (see Theorem \ref{t1}). Sec. \ref{s3} is devoted to the proof of Theorem \ref{t1} which is based on the spectral analysis of the linearized problem by using mentioned alternative definitions from \cite{Tcmp} of the Lopatinski condition and the uniform Lopatinski condition for 1-shocks. At last, in Appendix \ref{sec:appA}, for the reader's convenience, we describe the derivation of the Lopatinski determinant for symmetric hyperbolic systems with the 1-shock property.

\section{Linearized stability problem and main result}
\label{s2}

We consider a constant solution
$(U^+,U^-,\varphi) =(\widehat{U}^+,\widehat{U}^-,0)$ of systems \eqref{21} and the boundary conditions \eqref{bc} associated with the rectilinear shock wave $x_1=0$:
\[
\widehat{U}^{\pm}=(\hat{p}^{\pm},\hat{v}^{\pm},\widehat{F}_1^{\pm},\widehat{F}_2^{\pm}), \quad \hat{\rho}^{\pm}=\rho (\hat{p}^{\pm}) >0,\quad \hat{c}^{\pm}=1/\sqrt{\rho'(\hat{p}^{\pm})}>0,
\]
\[
\hat{v}^{\pm}=(\hat{v}_1^{\pm},\hat{v}_2^{\pm}),\quad
\widehat{F}_1^{\pm} =(\widehat{F}_{11}^{\pm},\widehat{F}_{21}^{\pm}), \quad \widehat{F}_2^{\pm} =(\widehat{F}_{12}^{\pm},\widehat{F}_{22}^{\pm}),
\]
where all the hat values are given constants. We assume that $\hat{v}_1^{\pm}>0$. In view of the third condition in \eqref{bc}, $\hat{v}_2^+=\hat{v}_2^-$ and we can choose a reference frame in which
$\hat{v}_2^+=\hat{v}_2^-=0$. The rest constants satisfy the relations
\begin{equation}
\begin{split}
& \frac{\hat{\rho}^+}{\hat{\rho}^-}=\frac{\hat{v}_1^-}{\hat{v}_1^+},\quad \frac{\hat{\rho}^+}{\hat{\rho}^-}\big\{ (\hat{v}_1^+)^2-\big((\widehat{F}_{11}^+)^2+(\widehat{F}_{12}^+)^2\big)\big\}[\hat{\rho}]=[\hat{p}],\\ & \big[\widehat{F}_{2j}\big]=0,\quad \big[\hat{\rho}\widehat{F}_{1j}\big]=0
\end{split}
\label{sbc}
\end{equation}
following from \eqref{bc}, where $j=1,2$, $[\hat{\rho}]=\hat{\rho}^+-\hat{\rho}^-$, etc.

We also assume that the constant solution satisfies the Lax conditions \eqref{25}:
\begin{equation}
M_->\frac{M}{\sqrt{M^2-M_1^2}},
\label{Mach-}
\end{equation}
\begin{equation}
M_1<M<M_*,
\label{Mach}
\end{equation}
where $M_-={\hat{v}_1^-}/{\hat{c}^-}$ and $M={\hat{v}_1^+}/{\hat{c}^+}$ are the upstream downstream  Mach numbers respectively,
\[
M_1=\sqrt{\mathcal{F}_{11}^2+\mathcal{F}_{12}^2},\quad M_*=\sqrt{1+\mathcal{F}_{11}^2+\mathcal{F}_{12}^2},
\]
and $\mathcal{F}_{ij}=\widehat{F}_{ij}^+/\hat{c}^+$ are the components of the unperturbed scaled deformation gradient $\mathcal{F}=(\mathcal{F}_{ij})_{i,j=1,2}$ behind of the shock. Note that \eqref{Mach-} follows from the first inequality in \eqref{25} and relations \eqref{sbc}.

As is known, for 1-shocks all the characteristics of the linearized system for the perturbation $\delta U^-$ ahead of the shock are incoming in the shock and
without loss of generality we may assume that $\delta U^-\equiv 0$. Following \cite{MTT20}, we write down the constant coefficient linearized  problem in a dimensionless form for the scaled perturbation $U=(p,v,F_1,F_2)$ behind of the shock wave resulting from the linearization of \eqref{bc}, \eqref{21} about the constant solution described above:
\begin{alignat}{3}
&Lp+{\rm div}\,v=0, &  \label{ls1}\\
& M^2Lv+\nabla p -(\mathcal{F}_1\cdot\nabla )F_1-(\mathcal{F}_2\cdot\nabla )F_2=0, & \label{ls2}\\
& LF_1- (\mathcal{F}_1\cdot\nabla )v=0,&  \label{ls3}\\
&  LF_2- (\mathcal{F}_2\cdot\nabla )v=0 &\qquad \mbox{for}\ x_1>0,\label{ls4}
\end{alignat}
\begin{alignat}{3}
& v_1 +d_0p -\frac{\ell_0}{M^2R}\,v_2=0, & \label{lbc1}\\
& a_0p+(1-R)\partial_{\star}\varphi =0,\qquad
v_2 +(1-R)\partial_2\varphi= 0, & \label{lbc2}\\
&  F_{11}+{\cal F}_{11}\,p-\frac{{\cal F}_{21}}{R}\,v_2=0, \qquad
 F_{12}+{\cal F}_{12}\,p-\frac{{\cal F}_{22}}{R}\,v_2=0,  & \label{lbc3}\\
& F_{21}-{\cal F}_{11}\,v_2=0, \qquad   F_{22}-{\cal F}_{12}\,v_2=0 & \qquad \mbox{on}\ x_1=0,\label{lbc4}
\end{alignat}
where
\[
L =\partial_t+\partial_1,\quad \mathcal{F}_1=(\mathcal{F}_{11},\mathcal{F}_{21}),\quad \mathcal{F}_2=(\mathcal{F}_{12},\mathcal{F}_{22}),\quad
d_0=\frac{M_*^2+M^2}{2M^2},\quad R=\frac{\hat{\rho}^+}{\hat{\rho}^-},
\]
\[
\ell_0={\cal F}_{11}{\cal F}_{21}+{\cal F}_{12}{\cal F}_{22},\quad  a_0=-\frac{\beta^2R}{2M^2},\quad
\beta =\sqrt{M_*^2-M^2}\quad (\mbox{cf.}\ \eqref{Mach}),\quad \partial_{\star}=\partial_t-\frac{\ell_0}{M^2}\,\partial_2,
\]
and we use the scaled values
\[
x'=\frac{x}{l},\quad t'=\frac{\hat{v}_1^+t}{l},\quad p=\frac{\delta p^+}{\hat{\rho}^+(\hat{c}^+)^2},\quad v=\frac{\delta v^+}{\hat{v}_1^+},\quad F_{ij}=\frac{\delta F_{ij}^+}{\hat{c}^+},\quad  \varphi =\frac{\delta\varphi}{l},
\]
for the original (unscaled) perturbation $\delta U^+=(\delta p^+,\delta v^+,\delta F_1^+,\delta F_2^+)$ behind the shock and the shock perturbation $\delta\varphi$, with $l$ being a typical length (the primes were dropped in \eqref{ls1}--\eqref{lbc4}). Our constant coefficient linearized problem is \eqref{ls1}--\eqref{lbc4} with the initial data
\begin{equation}
{U} (0,{x})={U}_0({x}),\quad {x}\in \mathbb{R}^2,\quad \varphi (0,{x}_2)=\varphi _0({x}_2),\quad {x}_2\in\mathbb{R}.\label{lindat}
\end{equation}


Let us recall the terminology. In our case, the rectilinear shock is called uniformly stable if the linearized problem \eqref{ls1}--\eqref{lindat} satisfies the uniform Lopatinski condition \cite{Kreiss}. This shock is called weakly stable if problem \eqref{ls1}--\eqref{lindat} satisfies the Lopatinski condition in a weak sense, i.e., problem \eqref{ls1}--\eqref{lindat} satisfies the Lopatinski condition \cite{Kreiss} but violates the uniform Lopatinski condition. We are now in a position to formulate the main result of this paper.

\begin{theorem}
Let a rectilinear shock wave in 2D compressible isentropic elastodynamics satis\-fies the Lax conditions \eqref{Mach-} and \eqref{Mach}. Then, this shock wave is uniformly stable if and only if
\begin{equation}
\begin{split}
\big(1+&\mathcal{F}_{11}^2+\mathcal{F}_{12}^2+M^2\big)\left(1+(\mathcal{F}:\mathcal{F})+(\det \mathcal{F})^2\right) \\ & -\left\{R(M^2-\mathcal{F}_{11}^2-\mathcal{F}_{12}^2)+\mathcal{F}_{21}^2+\mathcal{F}_{22}^2 \right\}\left(1+\mathcal{F}_{11}^2+\mathcal{F}_{12}^2\right)^2 \\ &
+\left({\cal F}_{11}{\cal F}_{21}+{\cal F}_{12}{\cal F}_{22}\right)^2\left\{2(1+\mathcal{F}_{11}^2+\mathcal{F}_{12}^2)-M^2\right\} \\
& \qquad >2M\left|{\cal F}_{11}{\cal F}_{21}+{\cal F}_{12}{\cal F}_{22}\right|\sqrt{\left(1+\mathcal{F}_{11}^2+\mathcal{F}_{12}^2-M^2\right)
 \left(1+(\mathcal{F}:\mathcal{F})+(\det \mathcal{F})^2\right)} \,,
\end{split}
\label{usc'}
\end{equation}
where $M$ is the downstream Mach number, $R$ measures the competition between downstream and upstream densities and $\mathcal{F}=(\mathcal{F}_{ij})_{i,j=1,2}$ is the scaled deformation gradient behind of the shock. If \eqref{usc'} is violated, then the rectilinear shock wave is weakly stable.
\label{t1}
\end{theorem}

An a priori estimate without loss of derivatives from the initial data was derived in \cite{MTT20} by constructing a so-called strictly dissipative 2-symmetrizer \cite{Tsiam}, provided that the stability condition \eqref{usc'} holds. Moreover, referring to \cite{Tsiam}, one can also write down an a priori estimate for the corresponding inhomogeneous problem, i.e., for problem \eqref{ls1}--\eqref{lindat} with a given source term $f(t,x)\in \mathbb{R}^7$ in the right-hand side of the interior equations \eqref{ls1}--\eqref{ls4} and a given source term $g(t,x_2)\in \mathbb{R}^7$ in the right-hand side of the boundary conditions \eqref{lbc1}--\eqref{lbc4}. This estimate reads
\begin{equation}
\begin{split}
\| U\|_{H^2([0,T]\times\mathbb{R}^2_+)}+ & \| U_{|x_1=0}\|_{H^2([0,T]\times\mathbb{R})} + \| \varphi\|_{H^3([0,T]\times\mathbb{R})}\\ &
\leq C\big\{\| U_0\|_{H^2(\mathbb{R}^2_+)} +\| \varphi_0\|_{L^2(\mathbb{R})} +\| f\|_{H^2([0,T]\times\mathbb{R}^2_+)}
+\| g\|_{H^2([0,T]\times\mathbb{R})}\big\},
\end{split}
\label{aprest'}
\end{equation}
where $\mathbb{R}^2_{+}=\{x_1>0,\ x_2\in\mathbb{R}\}$, the constant $C>0$ depends on $T$ and does not depend on the initial data and the source terms. Since estimate \eqref{aprest'} is an a priori estimate {\it without loss of derivatives} from the initial data and the source terms, the energy method in \cite{MTT20} can be considered as an indirect proof that condition \eqref{usc'} is sufficient for uniform stability.

By introducing the ``elastic'' Mach number $\widetilde{M}=\sqrt{M^2-M_1^2}\in (0,1)$, cf. \eqref{Mach}, and following \cite{MTT20}, we can equivalently rewrite \eqref{usc'} as
\begin{equation}\label{usc1}
\widetilde{M}^2(R-1)<1+ \frac{\mathcal{D}}{M_*^4},
\end{equation}
where
\[
\mathcal{D}=(M\sigma -|\ell_0|\beta -M_*^2\widetilde{M})(M\sigma -|\ell_0|\beta +M_*^2\widetilde{M}),
\]
with
\[
\sigma = \sqrt{M_*^2(1+M_2^2)-\ell_0^2}=\sqrt{M_*^2+M_2^2+(\det \mathcal{F})^2}\quad\mbox{and}\quad
M_2=\sqrt{\mathcal{F}_{21}^2+\mathcal{F}_{22}^2}.
\]

Setting formally $\mathcal{F}=0$ in \eqref{usc1}, we obtain the uniform stability condition
\begin{equation}
M^2(R-1)<1
\label{gas}
\end{equation}
found by Majda \cite{M1} (and written in our notations) for shock waves in isentropic gas dynamics. It was shown in \cite{MTT20} that $\mathcal{D}>0$. Comparing \eqref{usc1} and \eqref{gas}, this means that the elastic force plays stabilizing role. Moreover, it was proved in \cite{MTT20} that $\widetilde{M}^2(R-1)<1$ for compressive shock waves ($R>1$) with convex equations of state $p=p(\rho )$. This implies that these shock waves are uniformly stable because \eqref{usc1} holds for them. By the way, the rarefaction shock waves ($R<1$) are also always uniformly stable, cf. \eqref{usc1}.

For the particular deformations for which $\mathcal{F}_{12}=\mathcal{F}_{21}=0$ (pure stretching) or $\mathcal{F}_{11}=\mathcal{F}_{22}=0$, the stability condition \eqref{usc'} was shown in \cite{MTT20} by spectral analysis to be not only sufficient but also necessary for uniform stability. Actually, the spectral analysis in \cite{MTT20} for these two particular cases could be easily generalized to the case $\ell_0=0$ corresponding to a diagonal Cauchy--Green stress tensor $\mathcal{F}\mathcal{F}^\top$. Our goal now is to prove that \eqref{usc'} is necessary for uniform stability for general deformations. Moreover, we will show that violent instability never happens, that is maybe evident from the physical point of view.

\section{Proof of Theorem \eqref{t1}}
\label{s3}

Standard definitions of the Lopatinski condition and the uniform Lopatinski condition for initial boundary value problems for linear constant coefficient hyperbolic systems were given by Kreiss \cite{Kreiss}. For symmetric hyperbolic systems having the 1-shock property  one can use equivalent definitions. They were introduced in \cite{Tcmp} by using ideas of the normal modes analysis in \cite{GK} for MHD shock waves.
Following \cite{Tcmp}, we now give them for our problem \eqref{ls1}--\eqref{lindat}. For the reader's convenience, in Appendix \ref{sec:appA} we give the derivation of the Lopatinski determinant for symmetric hyperbolic systems with the 1-shock property (some points of the analysis are described in Appendix \ref{sec:appA} even in more detail than in \cite{Tcmp}).

We first rewrite system \eqref{ls1}--\eqref{ls4} in the matrix form
\begin{equation}\label{71}
\mathcal{A}_0\partial_tU+\mathcal{A}_1\partial_1U+\mathcal{A}_2\partial_2U=0\quad\mbox{for}\ x_1>0,
\end{equation}
where $\mathcal{A}_0={\rm diag}\,(1,M^2I_2,I_4)$,
\[
\mathcal{A}_1=\begin{pmatrix}
1 & e_1 & \underline{0} & \underline{0}  \\[7pt]
e_1^{\top}&M^2I_2 &  -\mathcal{F}_{11}I_2 & - \mathcal{F}_{12}I_2  \\[3pt]
\underline{0}^{\top} &- \mathcal{F}_{11}I_2 & I_2 & O_2 \\
\underline{0}^{\top} &- \mathcal{F}_{12}I_2 & O_2 & I_2
\end{pmatrix},\quad
\mathcal{A}_2=\begin{pmatrix}
0 & e_2 & \underline{0} & \underline{0}  \\[7pt]
e_2^{\top}&O_2 & - \mathcal{F}_{21}I_2 & - \mathcal{F}_{22}I_2  \\[3pt]
\underline{0}^{\top} &- \mathcal{F}_{21}I_2 & O_2 & O_2 \\
\underline{0}^{\top} &- \mathcal{F}_{22}I_2 & O_2 & O_2
\end{pmatrix}.
\]
By cross differentiation the perturbation $\varphi$ of the shock front can be excluded from the boundary conditions \eqref{lbc2}:
\begin{equation}\label{lbc2'}
\partial_{\star}v_2=a_0\partial_2p\quad \mbox{on}\ x_1=0.
\end{equation}
Then, the boundary conditions \eqref{lbc1}, \eqref{lbc3}, \eqref{lbc4}, \eqref{lbc2'} form the system
\begin{equation}
\mathfrak{B}_0\partial_tU+\mathfrak{B}_2\partial_2U+\mathfrak{B}_3U =0\quad \mbox{on}\ x_1=0, \label{72}
\end{equation}
where the matrices $\mathfrak{B}_{\alpha}$ ($\alpha =\overline{0,2}$) of order $6\times 7$ can be easily written down.

Following \cite{GK,Tcmp} (see also \eqref{A9}, \eqref{A10} in Appendix \ref{sec:appA}), we write down the following linear algebraic system for a vector $X$ associated with problem \eqref{71}, \eqref{72}:
\begin{align}
& (s\mathcal{A}_0+\lambda^+ \mathcal{A}_1+{i}\omega \mathcal{A}_2)X =0\,,\label{3.5}
\\
& \big(\mathcal{A}_1\widetilde{U}_0\big)\cdot X=0,
\label{3.6}
\end{align}
where $s=\eta +{i}\xi$, $\eta >0$, $(\xi ,{\omega})\in \mathbb{R}^2$ ($s$ and $\omega$ are, in fact, the Laplace and Fourier variables respectively), the vector $\widetilde{U}_0=\widetilde{U}_0(s,\omega )$ satisfies
\begin{equation}
(s\mathfrak{B}_0+i\omega\mathfrak{B}_2+\mathfrak{B}_3)\widetilde{ U}_0=0,\label{3.2}
\end{equation}
and $\lambda^+=\lambda^+ (\eta ,\xi ,\omega )$ is a simple root $\lambda$ of the dispersion relation
\begin{equation}
\det (s\mathcal{A}_0+\lambda \mathcal{A}_1+{i}\omega \mathcal{A}_2)=0 \label{3.4}
\end{equation}
lying strictly in the open right-half complex plane ($\Real\lambda >0$). Note that since our shock waves are 1-shocks, the boundary matrix $\mathcal{A}_1$ has only one negative eigenvalue. Then, in view of Hersh's lemma \cite{Hersh}, for all $\eta >0$ and $(\xi ,{\omega})\in \mathbb{R}^2$ equation \eqref{3.4} has a {\it unique} solution $\lambda =\lambda^+$ with $\Real\lambda >0$.

Since $\lambda^+$ is a simple root, we can choose six linearly independent equations of system \eqref{3.5}. Adding them to equation \eqref{3.6}, we obtain for the vector $X$ the linear system
\begin{equation}
\mathfrak{L}X =0
\label{X}
\end{equation}
whose determinant is, in fact, the Lopatinski determinant (see \cite{Tcmp} and Appendix \ref{sec:appA}). We are now ready to give, as in \cite{Tcmp}, the definitions of the Lopatinski condition and the uniform Lopatinski condition for problem \eqref{71}, \eqref{72}.

\begin{definition}
Problem \eqref{71}, \eqref{72} satisfies the Lopatinski condition if the Lo\-pa\-tinski determinant $\det \mathfrak{L} (\eta
,\xi ,\omega ,\lambda^+ ) \neq 0$ for all $\eta >0$ and  $(\xi , \omega )\in
\mathbb{R}^2$.
\label{d3.2}
\end{definition}

\begin{definition}
Problem \eqref{71}, \eqref{72} satisfies the uniform Lopatinski condition if the Lopatinski determinant $\det \mathfrak{L} (\eta
,\xi ,\omega ,\lambda^+ ) \neq 0$ for all $\eta \geq 0$ and  $(\xi , \omega )\in
\mathbb{R}^2$ (with $\eta^2+\xi^2+{\omega}^2\neq 0$), where $\lambda^+ (0,\xi ,\omega )= \lim\limits_{\eta \rightarrow +0}
\lambda^+ (\eta ,\xi ,\omega )  $.
\label{d3.3}
\end{definition}

We write down the dispersion relation \eqref{3.4}, which is  a polynomial equation for finding the unique $\lambda=\lambda^+$:
\begin{equation}
\Omega^3\left( M^2\Omega^2-\sigma_1^2-\sigma_2^2\right)\left( M^2\Omega^2-\sigma_1^2-\sigma_2^2-\lambda^2+\omega^2\right)=0,
\label{3.4'}
\end{equation}
with $\Omega =s+\lambda$, $\sigma_1=\mathcal{F}_{11}\lambda+i\omega\mathcal{F}_{21}$ and $\sigma_2=\mathcal{F}_{12}\lambda+i\omega\mathcal{F}_{22}$. As was shown in \cite{MTT20},  $\lambda^+$ is one of the two roots of the equation
\begin{equation}
 M^2\Omega^2-\sigma_1^2-\sigma_2^2-\lambda^2+\omega^2=0
\label{3.4''}
\end{equation}
whose left-hand side is the last multiplier in the left-hand side of \eqref{3.4'}.

Omitting straightforward calculations, we find the vector $\widetilde{U}_0$ in \eqref{3.2}, which is determined up to a nonzero factor:
\begin{multline*}
\widetilde{U}_0=\bigg(s-\frac{i\ell_0}{M^2}\omega  \,,\;-d_0s+\frac{i\ell_0}{M^2}\omega\,,\;ia_0\omega\,,\;-\mathcal{F}_{11}s+i\Big(\frac{a_0\mathcal{F}_{21}}{R}+
\frac{\ell_0\mathcal{F}_{11}}{M^2}\Big)\omega\,, \\
 ia_0\mathcal{F}_{11}\omega\,,\;i\Big(\frac{a_0\mathcal{F}_{22}}{R}+\frac{\ell_0\mathcal{F}_{12}}{M^2}\Big)\omega\,,\; ia_0\mathcal{F}_{12}\omega\bigg).
\end{multline*}
We can then calculate the vector $\mathcal{A}_1\widetilde{U}_0$ appearing in \eqref{3.6}:
\[
\mathcal{A}_1\widetilde{U}_0=-\frac{\beta^2}{2M^2} \big( s\,,\;i\ell_0\omega-sM^2\,,\;iR(M-M_1^2)\omega\,,\;i\mathcal{F}_{21}\omega-\mathcal{F}_{11}s\,,\;0\,,\;
i\mathcal{F}_{22}\omega-\mathcal{F}_{12}s\,,\;0\big).
\]

Replacing the first line of the matrix $s\mathcal{A}_0+\lambda^+ \mathcal{A}_1+{i}\omega \mathcal{A}_2$ with the vector $\mathcal{A}_1\widetilde{U}_0$ considered  as a row-vector, we get the matrix $\mathfrak{L}$ (see \eqref{3.5}, \eqref{3.6}). Omitting calculations, we obtain
\[
\begin{split}
\det \mathfrak{L} =\frac{\beta^2\Omega^2(\omega^2-\lambda^2)}{2M^2}\Big\{  (\lambda^2-\omega^2)s & +(M^2s-i\ell_0\omega )\Omega\lambda  +M_{1}^2\lambda^2s \\ & + M_{2}^2\omega^2\lambda +i\ell_0\omega\lambda (s-\lambda )+R(M^2-M_{1}^2)\omega^2\Omega \Big\},
\end{split}
\]
where
$\lambda =\lambda^+$, i.e., $\lambda$ should be the solution of equation \eqref{3.4''} with the property $\Real\lambda >0$ for $\eta >0$. Since
$\Omega^2(\omega^2-\lambda^2)\neq 0$ for $\lambda =\lambda^+$, the equality $\det \mathfrak{L}=0$ is equivalent to
\begin{equation}
\begin{split}
(\lambda^2-\omega^2)s  +  (M^2s-i\ell_0\omega )\Omega\lambda  +M_{1}^2\lambda^2s & + M_{2}^2\omega^2\lambda \\ & +i\ell_0\omega\lambda (s-\lambda )+R(M^2-M_{1}^2)\omega^2\Omega =0.
\end{split}
\label{3.4'''}
\end{equation}
The test of the (uniform) Lopatinski condition is thus reduced to the study of solutions $(s,\lambda )$ of  system \eqref{3.4''}, \eqref{3.4'''} for all real $\omega$.

In view of \eqref{Mach} and the fact that
\[
M^2\sigma^2-\ell_0^2\beta^2=M_*^2\left(M_2^2(M^2-M_1^2)+M^2+(\det \mathcal{F})^2\right)>0,
\]
the following three parameters are strictly positive:
\[
K=R(M^2-M_1^2)+M_2^2 >0,\quad K_1=\frac{(M\sigma -|\ell_0|\beta )^2}{M_*^4}>0,\quad K_2=1+M_2^2 >0.
\]
For the particular cases $\mathcal{F}_{12}=\mathcal{F}_{21}=0$ and $\mathcal{F}_{11}=\mathcal{F}_{22}=0$ (for which $\ell_0=0$) studied in \cite{MTT20} by spectral analysis, the domains of uniform and weak stability can be described in terms of the above parameters $K$, $K_1$ and $K_2$. But, as we will see, the same is true for general deformations. We can already check that the stability condition \eqref{usc'} is equivalently reformulated as
\begin{equation}\label{usc^}
 K < K_1 +K_2.
\end{equation}
It was proved in \cite{MTT20} by the energy method that \eqref{usc^} is sufficient for uniform stability. For simplifying the arguments of our spectral analysis, we will below use this knowledge. Let us assume that
\begin{equation}\label{KK}
K\geq K_1+K_2.
\end{equation}

The dispersion relation \eqref{3.4''} can be rewritten as
\begin{equation}\label{83'}
M^2\Omega^2-M_*^2\lambda^2 +K_2\omega^2=2i\ell_0\lambda\omega
\end{equation}
whereas equality \eqref{3.4'''} reads
\begin{equation}\label{84}
\Omega (M^2\lambda s +K\omega^2) +(M_*^2\lambda^2  -K_2\omega^2)s=2i\ell_0\lambda^2\omega .
\end{equation}
It follows from \eqref{83'} that $M_*^2\lambda^2-K_2\omega^2 =M^2\Omega^2 -2i\ell_0\lambda\omega$. Substituting this into \eqref{84}, we obtain
\begin{equation}
\Omega (M^2\Omega^2-M^2\lambda^2+K\omega^2-2i\ell_0\lambda\omega )=0.
\label{0.}
\end{equation}
Since $\Omega\neq 0$ for $\lambda =\lambda^+$ and $\eta >0$, \eqref{0.} is reduced to
\[
M^2\Omega^2-M^2\lambda^2+K\omega^2=2i\ell_0\lambda\omega .
\]

We can thus consider the system
\begin{align}
 M^2\Omega^2-M_*^2\lambda^2 +K_2\omega^2 & =2i\ell_0\lambda\omega, \label{83=}\\
 M^2\Omega^2-M^2\lambda^2+K\omega^2 & =2i\ell_0\lambda\omega \label{86=}
\end{align}
instead of system \eqref{3.4''}, \eqref{3.4'''}. Moreover, we can simplify our arguments below by avoiding a separate consideration of the cases $\ell_0 \geq 0$ and $\ell_0<0$. Namely, instead of system \eqref{83=}, \eqref{86=} we may consider the system
\begin{align}
 M^2\Omega^2-M_*^2\lambda^2 +K_2\omega^2 & =2i|\ell_0|\lambda\omega, \label{83}\\
 M^2\Omega^2-M^2\lambda^2+K\omega^2 & =2i|\ell_0|\lambda\omega .\label{86}
\end{align}
Indeed, these systems coincide if $\ell_0\geq 0$.  On the other hand, if $\ell_0< 0$, then in \eqref{83=}, \eqref{86=} we make the change $\widetilde{\omega}=-\omega \in \mathbb{R}$.  After dropping tildes we again get system \eqref{83}, \eqref{86}.

The dispersion relation \eqref{83} has the two roots
\[
\lambda_k=\frac{1}{\beta^2}\left( M^2s-i|\ell_0|\omega +(-1)^k\sqrt{M^2M_*^2s^2-2i|\ell_0|M^2s\omega+(K_2\beta^2-\ell_0^2)\omega^2}\right),\quad k=1,2.
\]
In view of \eqref{Mach}, we have:
\[
\Real\lambda_1|_{\omega =0} =\frac{M}{\beta^2}(M-M_*)\eta <0,\quad \Real\lambda_2|_{\omega =0} =\frac{M}{\beta^2}(M+M_*)\eta >0\quad\mbox{for}\quad \eta>0.
\]
Hence, by virtue of Hersh's lemma \cite{Hersh}, $\Real\lambda_2 >0$ for all $\omega\in\mathbb{R}$, i.e.,
\begin{equation}
\lambda^\pm=\frac{1}{\beta^2}\left( M^2s-i|\ell_0|\omega \pm\sqrt{M^2M_*^2s^2-2i|\ell_0|M^2s\omega+(K_2\beta^2-\ell_0^2)\omega^2}\right),
\label{l+}
\end{equation}
where $\lambda^-:=\lambda_1$. It follows from \eqref{83}, \eqref{86} that
\begin{equation}\label{l2}
\lambda^2=\frac{(K_2-K)\omega^2}{\beta^2}.
\end{equation}
By virtue of \eqref{KK}, this implies $\Real\lambda =0$. Since $\Real\lambda^+ >0$ for $\eta >0$, we necessarily have that $\eta =0$, i.e., shock waves cannot be violently unstable and in the parameter domain \eqref{KK}  they are, at least, weakly stable.

For $\omega =0$, \eqref{l+} and \eqref{l2} yield $s =0$. Hence, the uniform Lopatinski condition holds for the 1D case because the solution $s=0$ is prohibited by the requirement $\eta^2+\xi^2+\omega^2 \neq 0$ (see Definition \ref{d3.3}). That is, we may assume that $\omega \neq 0$. Since the left-hand sides in \eqref{83} and \eqref{86} are homogeneous functions of $s$, $\lambda$ and $\omega$, without loss of generality, from now on we will suppose that $\omega =1$.

For $\eta =0$,
\begin{equation}\label{iml}
\Imag\lambda^\pm =\delta^\pm =\frac{1}{\beta^2}\left( M^2\xi-|\ell_0| \pm\sqrt{M^2M_*^2\xi^2-2|\ell_0|M^2\xi+\ell_0^2-K_2\beta^2}\right).
\end{equation}
In the parameter domain of uniform stability, system \eqref{83}, \eqref{86} has no roots $(s,\lambda )= (i\xi ,i\delta^+)$. Because of the continuous dependence of this system on $s$ and $\lambda$, the passage to weak stability may happen only thanks to the merging of $\delta^+$ and $\delta^-$ at a point of transition $\xi =\xi_*$. It is clear that for $\xi=\xi_*$ the square root in \eqref{iml} should vanish:
\begin{equation}
M^2M_*^2\xi^2_*-2|\ell_0|M^2\xi_*+\ell_0^2-K_2\beta^2 =0.
\label{xi2}
\end{equation}
From \eqref{xi2} we find
\begin{equation}\label{xi*}
\xi_*=\xi_*^\pm =\frac{M|\ell_0|\pm\beta\sigma}{MM_*^2},
\end{equation}
and we have
\begin{equation}\label{d*}
\delta^+_{|\xi =\xi_*}=\delta^-_{|\xi =\xi_*}=  \delta_*=\frac{M^2\xi_*-|\ell_0|}{\beta^2}\qquad (\mbox{for}\ \xi_*=\xi_*^\pm ).
\end{equation}
Since $\delta^\pm$ is real, the elementary analysis of the quadratic function in \eqref{xi2} shows that
\begin{equation}
\xi \geq \xi_*^+\quad\mbox{or}\quad \xi\leq\xi_*^-.
\label{xilim}
\end{equation}

Substituting \eqref{xi*} into \eqref{d*}, we obtain
\begin{equation}
\delta_* =\frac{\sqrt{K_1}}{\beta} \quad\mbox{for} \ \xi_*=\xi_*^+\qquad\mbox{and}\qquad
\delta_* =-\frac{\sqrt{K_3}}{\beta} \quad\mbox{for} \ \xi_*=\xi_*^-,
\label{delt}
\end{equation}
where
\[
K_3=\frac{(M\sigma +|\ell_0|\beta )^2}{M_*^4}>K_1>0.
\]

From \eqref{l2} and \eqref{delt} we find (for $K>K_2$, cf. \eqref{KK}) two possible transitions to weak stability:
\[
K=K_1+K_2\quad \mbox{and}\quad K=K_3+K_2.
\]
Here the point $K=K_1+K_2$ (cf. \eqref{usc^} and \eqref{KK}) on the $K$--axis corresponds to $\xi_*=\xi_*^+$ and lies to the left from the point $K=K_3+K_1$ corresponding to $\xi_*=\xi_*^-$. Recall that at the left from the point $K=K_1+K_2$ on the $K$--axis we have the domain of uniform stability, see \eqref{usc^} (this was proved in \cite{MTT20} by the energy method). Recall also that for the domain of weak stability system \eqref{83}, \eqref{86} necessarily has a root  $(s,\lambda )= (i\xi ,i\delta^+)$. It follows from \eqref{iml} and \eqref{xilim} that
\[
\delta^+ \geq \frac{M^2\xi-|\ell_0|}{\beta^2} \geq \frac{M^2\xi_*^+ -|\ell_0|}{\beta^2}=\frac{\sqrt{K_1}}{\beta},
\]
or $(\delta^+)^2\geq K_1/\beta^2$. In view of \eqref{l2}, the last inequality is written as $K\geq K_1+K_2$, cf. \eqref{KK}.  Note that, by virtue of \eqref{iml} and \eqref{xilim}, the point $K=K_3+K_2$ on the $K$--axis just corresponds to the emergence of the root  $(s,\lambda^- )= (i\xi ,i\delta^-)$  for $K\geq K_3+K_1$ whereas system \eqref{83}, \eqref{86} has the ``weakly stable'' root  $(s,\lambda^+ )= (i\xi ,i\delta^+)$ for all $K\geq K_1+K_2$ (in particular, for  $K\geq K_3+K_1$). That is, inequality \eqref{KK} describes the parameter domain of weak stability. This completes the proof of Theorem \ref{t1}.

\section*{\normalsize Acknowledgements}

This research was supported by the Russian Science Foundation under grant No. 20-11-20036.

\appendix
\titleformat{\section}{\large\bfseries}{Appendix \thesection}{1em}{}

\section{Lopatinski determinant for symmetric hyperbolic systems with the 1-shock property} \label{sec:appA}

Let us consider an abstract $n$-dimensional counterpart of problem \eqref{71}, \eqref{72}:
\begin{align}
\mathcal{A}_0\partial_tU+\sum_{j=1}^n\mathcal{A}_j\partial_jU=0 & \quad\mbox{for}\ x_1>0,
\label{A1} \\
\mathfrak{B}_0\partial_tU+\sum_{k=2}^n\mathfrak{B}_k\partial_kU+\mathfrak{B}_{n+1}U =0 &\quad \mbox{on}\ x_1=0,
\label{A2}
\end{align}
where $U=U(t,x)\in \mathbb{R}^m$, $x\in\mathbb{R}^n$, $\mathcal{A}_{\alpha}\in \mathbb{M}(m,m)$ ($\alpha =\overline{0,n}$) are symmetric matrices, $\mathcal{A}_0>0$ and $\mathfrak{B}_{\alpha}\in \mathbb{M}(m-1,m)$ ($\alpha =\overline{0,n+1}$). Moreover, we assume that the boundary $x_1=0$ is not characteristic ($\det \mathcal{A}_1\neq 0$) and the matrix $\mathcal{A}_1$ has one negative and $m-1$ positive eigenvalues.

We construct an Hadamard-type ill-posedness example
\begin{equation}
U_k=e^{-\sqrt{k}+ k(st+i(\omega'\cdot x'))}\widetilde{U}(kx_1)\qquad (k\in\mathbb{N})
\label{A3}
\end{equation}
for problem \eqref{A1}, \eqref{A2}, where
\[
s=\eta +{i}\xi ,\quad \eta >0,\quad (\xi ,{\omega}')\in \mathbb{R}^n,\quad \omega'=(\omega_2,\ldots ,\omega_n),\quad  x'=(x_2,\ldots ,x_n).
\]
The substitution of \eqref{A3} into \eqref{A1} and \eqref{A2} gives the boundary value problem
\begin{align}
& \frac{{\rm d}\widetilde{U}}{{\rm d} x_1}=\mathfrak{A}(s,{\omega}' )\widetilde{ U}, \quad
 x_1>0,\label{A4}\\
& \Big(s\mathfrak{B}_0+i\sum_{k=2}^n\omega_k\mathfrak{B}_k+\mathfrak{B}_{n+1}\Big)\widetilde{ U}_0=0\label{A5}
\end{align}
for $\widetilde{U}(x_1)$, where
\[
\mathfrak{A}=\mathfrak{A}(s,\omega' )=-\mathcal{A}_1^{-1}\mathcal{A},\quad \mathcal{A}= \mathcal{A}(s,\omega' )=s\mathcal{A}_0+{i}\sum_{k=2}^{n}\omega_k \mathcal{A}_k,\quad
\widetilde{ U}_0=\widetilde{ U}(0).
\]

Applying the Laplace transform to \eqref{A4}, we get $\lambda V-\widetilde{U}_0=\mathfrak{A}V$ or $(\mathcal{A}+\lambda\mathcal{A}_1)V=\mathcal{A}_1\widetilde{U}_0$ for
\[
V=V(\lambda )=\int_{0}^{+\infty}e^{-\lambda x_1}\widetilde{U}(x_1)dx_1.
\]
Then the inverse Laplace transform gives
\[
\widetilde{U}(x_1)
=\frac{1}{2\pi i}\oint\limits_{C} e^{\lambda x_1}(\mathcal{A}(s,\omega' )+\lambda \mathcal{A}_1)^{-1}
\mathcal{A}_1\widetilde{U}_0\,d\lambda ,
\]
where $C$ is a contour large enough to enclose all the singularities  of the integrand $F(\lambda )=e^{\lambda x_1}f(\lambda )$, with
$f(\lambda )=(\mathcal{A}(s,\omega' )+\lambda \mathcal{A}_1)^{-1}\mathcal{A}_1\widetilde{U}_0$. These singularities are the eigenvalues $\lambda$ of $\mathfrak{A}$, which satisfy
\begin{equation}
\det \Big(s\mathcal{A}_0+\lambda\mathcal{A}_1+{i}\sum_{k=2}^{n}\omega_k \mathcal{A}_k\Big)=0.
\label{A7}
\end{equation}

Then,
\[
\widetilde{U}(x_1)=\sum_{j=1}^l\res_{\lambda=\lambda_j} F(\lambda),
\]
where $l\leq m$ is the number of the eigenvalues $\lambda =\lambda_j$ without counting multiplicity.
In view of Hersh's lemma \cite{Hersh}, for all $\eta >0$ and $(\xi ,{\omega}')\in \mathbb{R}^n$ equation \eqref{A7} has a {\it unique} solution $\lambda =\lambda^+=\lambda^+ (\eta ,\xi ,\omega )$ lying strictly in the open right-half complex plane ($\Real\lambda >0$). For $\lambda = \lambda^+$ we have $e^{\lambda x_1}\rightarrow +\infty$ as $x_1\rightarrow +\infty$. We have a bounded $\widetilde{U}$ in \eqref{A3} and thus construct an Hadamard-type ill-posedness example if and only if the residue at this simple $\lambda^+$ is zero. Since $\lambda^+$ is a simple pole of the integrand $F(\lambda )$, then
$\res_{\lambda=\lambda^+}F(\lambda) =e^{\lambda^+ x_1}\res_{\lambda=\lambda^+}f(\lambda)$. We have $\res_{\lambda=\lambda^+}f(\lambda)=0$ if and only if $f (\lambda )$ is bounded at $\lambda =\lambda^+$. This is true if and only if there exists such a bounded vector $Y$ that
\begin{equation}
(\mathcal{A}(s,\omega' )+\lambda \mathcal{A}_1) Y = \mathcal{A}_1\widetilde{U}_0,
\label{A8}
\end{equation}
where $\lambda =\lambda^+$ is our unique root of \eqref{A7} with $\Real\lambda >0$.

A left annihilator $X$ of the matrix $\mathcal{A}(s,\omega' )+\lambda \mathcal{A}_1$ should also annihilate the right-hand side in \eqref{A8}: $X^{\top}\mathcal{A}_1\widetilde{U}_0=0$. Since the matrices are symmetric, we finally have
\begin{align}
& \Big(s\mathcal{A}_0+\lambda \mathcal{A}_1+{i}\sum_{k=2}^{n}\omega_k \mathcal{A}_k\Big)X =0,\label{A9}\\
& \big(\mathcal{A}_1\widetilde{U}_0\big)\cdot X=0,\label{A10}
\end{align}
where $\lambda =\lambda^+$ is the unique root of \eqref{A7} with $\Real\lambda >0$. Since $\lambda^+$ is a simple eigenvalue, we can choose $m-1$ linearly independent equations of system \eqref{A9}. Adding them to equation \eqref{A10}, we obtain for the vector $X$ a linear system
$\mathfrak{L}X =0$ whose determinant is the {\it Lopatinski determinant}. Indeed, problem \eqref{A1}, \eqref{A2} is ill-posed (we have so-called violent instability) if and only if the equation $\det \mathfrak{L} (s ,\omega ,\lambda^+ ) =0$ has a root $s$ with $\Real s=\eta >0$ for some $\omega'\in \mathbb{R}^{n-1}$.

\end{document}